\documentclass{article}
\usepackage{amsthm}
\usepackage{amssymb}
\usepackage{amsfonts}
\usepackage{latexsym}
\usepackage{amsmath}
\usepackage{indentfirst}

\newtheorem{thm}{Theorem}[section]
\newtheorem{lem}[thm]{Lemma}
\newtheorem{dfn}[thm]{Definition}
\theoremstyle{remark}
\newtheorem*{rmk}{Remark}
\newtheorem*{ack}{Acknowledgement}
\usepackage[top=30truemm,bottom=25truemm,left=25truemm,right=25truemm]{geometry}

\title{Arithmetic formulas for the Fourier coefficients of Hauptmoduln of level 2, 3, and 5}
\author{Toshiki Matsusaka and Ryotaro Osanai}
\date{}

\begin{document}
\maketitle

\begin{abstract}
We give arithmetic formulas for the coefficients of Hauptmoduln of higher levels as analogues of Kaneko's formula for the elliptic modular $j$-invariant. We also obtain their asymptotic formulas by employing Murty-Sampath's method.
\end{abstract}

\section{Introduction}

For the elliptic modular function $j(\tau)$, let  $\textbf{t}_m(d)$ be the modular trace function (the precise definition will be given later) and $c_n$ ($n \geq 1$) the $n$th Fourier coefficient of $j(\tau)$, that is, $j(\tau) = q^{-1} + 744 + \sum_{n=1}^{\infty} c_n q^n$. Zagier \cite{Zag02} studied the traces of singular moduli and showed that the generating function of $\textbf{t}_m(d)$ is a meromorphic modular form of weight $3/2$ on the right group for each $m$. Multiplying it by theta function and observing the modular forms of weight 2, Kaneko \cite{Kan95} gave the following arithmetic formula for $c_n$ experimentally, and showed it.
\begin{eqnarray*}
c_n &=& \frac{1}{n} \biggl\{ \sum_{r \in \mathbb{Z}} \textbf{t}_1(n - r^2) + \sum_{\substack{r \geq 1,\ odd}} ((-1)^n \textbf{t}_1(4n - r^2) - \textbf{t}_1(16n - r^2)) \biggr\} \\
&=& \frac{1}{2n} \sum_{r \in \mathbb{Z}} \textbf{t}_2(4n - r^2) .
\end{eqnarray*}
On the other hand, by using the circle method, Petersson \cite{Pet32} and later Rademacher \cite{Rad} independently derived the asymptotic formula for $c_n$:
\begin{eqnarray*}
c_n \sim \frac{e^{4 \pi \sqrt{n}}}{\sqrt{2} n^{3/4}}\ as\ n \to \infty.
\end{eqnarray*}
The circle method is introduced by Hardy and Ramanujan \cite{HR18} to prove the asymptotic formula for the partition function
\begin{eqnarray*}
p(n) \sim \frac{e^{\pi \sqrt{2n/3}}}{4\sqrt{3}n}\ as\ n \to \infty,
\end{eqnarray*}
where $p(n)$ is defined by $\sum_{n=0}^{\infty} p(n)q^n = \prod_{n=1}^{\infty}(1-q^n)^{-1}$. In 2013, Bruinier and Ono \cite{BO13} considered certain traces of singular moduli for weak Maass forms and derived the algebraic formula for $p(n)$. Combining this formula with Laplace's method, Dewar and Murty \cite{DM13, DM132} proved the asymptotic formulas for $p(n)$ and $c_n$ without the circle method. More recently, Murty and Sampath \cite{Sam15} derived the asymptotic formula for $c_n$ from Kaneko's arithmetic formula with Laplace's method.\\

In this article, we generalize these formulas to Hauptmoduln (defined in section 2) for the congruence subgroups $\Gamma_0(p)$ and $\Gamma_0^*(p)$ (the extension of $\Gamma_0(p)$ by the Atkin-Lehner involution) with $p = 2, 3,$ and $5$. \\

Let $j_p(\tau)$ and $j_p^*(\tau)$ be the corresponding Hauptmoduln for $\Gamma_0(p)$ and $\Gamma_0^*(p)$, respectively. Ohta \cite{Ohta09} gave the arithmetic formulas for the Fourier coefficients of $j_2(\tau)$ and $j_2^*(\tau)$, and a part of those of $j_3(\tau)$. She also treated the cases of $j_4(\tau)$ and $j_4^*(\tau)$. Let $c_n^{(p)}$ and $c_n^{(p*)}$ be the $n$th Fourier coefficients of $j_p(\tau)$ and $j_p^*(\tau)$, respectively. We express these coefficients in terms of the modular trace functions $\textbf{t}_m^{(p*)}(d)$.

\begin{thm}
\label{main1}
For any $n \geq 1$, we have
\begin{eqnarray*}
c_n^{(2)} &=& \frac{1}{2n} \times \left\{ \begin{array}{ll}
 - \sum_{r \equiv 0 (2)} \textbf{t}_2^{(2*)}(4n - r^2) + 24\sigma_1^{(2)}(n)  & (n \equiv 0 \bmod 2), \\
  \sum_{r \in \mathbb{Z}} \textbf{t}_2^{(2*)}(4n - r^2) + 24\sigma_1(n)  & (n \not\equiv 0 \bmod 2), \\ \end{array} \right. \\
c_n^{(3)} &=& \frac{1}{2n} \times \left\{ \begin{array}{ll}
 - \sum_{r \equiv 0 (3)} \textbf{t}_2^{(3*)}(4n - r^2) + 36\sigma_1^{(3)}(n)  & (n \equiv 0 \bmod 3), \\
 \sum_{r \in \mathbb{Z}} \textbf{t}_2^{(3*)}(4n - r^2) + 36\sigma_1(n) & (n \not\equiv 0 \bmod 3), \\ \end{array} \right. \\
c_n^{(5)} &=& \frac{1}{2n} \times \left\{ \begin{array}{ll}
 -\sum_{r \equiv 0 (5)} \textbf{t}_2^{(5*)}(4n - r^2) + 18\sigma_1^{(5)}(n)  & (n \equiv 0 \bmod 5), \\
\sum_{r \in \mathbb{Z}} \textbf{t}_2^{(5*)}(4n - r^2) +18\sigma_1(n) & (n \not\equiv 0 \bmod 5), \\ \end{array} \right. \\
c_n^{(p*)} &=& c_n^{(p)} - p c_{pn}^{(p)} \ \  (p = 2, 3, 5)\\
\end{eqnarray*}
where $\sigma_1(n) = \sum_{d | n}d$, and $\sigma_1^{(p)}(n) = \sum_{\substack{d | n \\ p \nmid d}}d$.
\end{thm}

\begin{rmk}
These formulas are different from those in Ohta \cite{Ohta09}. In \cite{Ohta09}, the definition of $\textbf{t}_m^{(p)}(d)$ was mixed with that of $\textbf{t}_m^{(p*)}(d)$, and used the values of $\textbf{t}_m^{(p)}(d)$ instead of $\textbf{t}_m^{(p*)}(d)$.
\end{rmk}
Combining these formulas with Laplace's method as in \cite{Sam15}, we obtain the asymptotic formulas of $c_n^{(p)}$.

\begin{thm}
\label{main2}
We have
\begin{eqnarray*}
c_n^{(2)} &\sim& \frac{e^{2\pi \sqrt{n}}}{2n^{3/4}} \times \left\{ \begin{array}{ll}
 -1 & (n \equiv 0 \bmod 2), \\ 1 & (n \equiv 1 \bmod 2), \\
\end{array} \right.\\
c_n^{(3)} &\sim& \frac{e^{4\pi \sqrt{n}/3}}{\sqrt{6}n^{3/4}} \times \left\{ \begin{array}{ll}
 -1 & (n \equiv 0, 2 \bmod 3), \\ 2 & (n \equiv 1\ \ \ \bmod 3), \\
\end{array} \right.\\
c_n^{(5)} &\sim& \frac{e^{4\pi \sqrt{n}/5}}{\sqrt{10}n^{3/4}} \times \left\{ \begin{array}{ll}
 -1 & (n \equiv 0 \bmod 5), \\ (3 + \sqrt{5})/2 & (n \equiv 1 \bmod 5), \\ -1 + \sqrt{5} & (n \equiv 2 \bmod 5), \\ -1 - \sqrt{5} & (n \equiv 3 \bmod 5), \\ (3 - \sqrt{5})/2 & (n \equiv 4 \bmod 5)\\
\end{array} \right.
\end{eqnarray*}\\
as n $\to \infty$.
\end{thm}

\section{Preliminaries}
In this section, we shall define the Hauptmoduln and the modular trace functions.

\begin{dfn}
Let $\Gamma$ be a congruence subgroup of $SL_2(\mathbb{R})/{\pm I}$ containing $(\begin{smallmatrix}1 & 1 \\0 & 1 \end{smallmatrix})$. If the genus of $\Gamma$ is equal to 0, there is a unique modular function f of weight 0 satisfying the following conditions. We call this f the Hauptmodul with respect to $\Gamma$.\\
$(1)$ f is holomorphic in the upper half plane $\mathfrak{H}$,\\
$(2)$ f has a Fourier expansion of the form $f(\tau) = q^{-1} + \sum_{n = 1}^{\infty} H_n q^n \ (q := e^{2\pi i \tau})$,\\
$(3)$ f is holomorphic at cusps of $\Gamma$ except i$\infty$.
\end{dfn}

For $\Gamma_0(p) := \{(\begin{smallmatrix}a & b \\c & d \end{smallmatrix}) \in \mathrm{PSL}_2(\mathbb{Z})\ |\ c \equiv 0 \pmod{p} \}$ and $\Gamma_0^*(p) := \Gamma_0(p) \cup \Gamma_0(p)(\begin{smallmatrix}0 & -1/\sqrt{p} \\ \sqrt{p} & 0 \end{smallmatrix})$ ($p$ = 2, 3, 5), the corresponding Hauptmoduln $j_p(\tau)$ and $j_p^*(\tau)$ can be described by means of the Dedekind $\eta$-function $\eta(\tau):= q^{1/24}\prod_{n=1}^{\infty}(1-q^n)$;
\begin{eqnarray*}
j_2(\tau) &=& \biggl( \frac{\eta(\tau)}{\eta(2\tau)} \biggr)^{24} + 24 = \frac{1}{q} + 276q - 2048q^2 + 11202q^3 + \cdots,\\
j_2^*(\tau) &=& \biggl( \frac{\eta(\tau)}{\eta(2\tau)} \biggr)^{24} + 24 + 2^{12} \biggl( \frac{\eta(2\tau)}{\eta(\tau)} \biggr)^{24} = \frac{1}{q} + 4372q + 96256q^2 + 1240002q^3 + \cdots,\\
j_3(\tau) &=& \biggl( \frac{\eta(\tau)}{\eta(3\tau)} \biggr)^{12} + 12 = \frac{1}{q} + 54q - 76q^2 - 243q^3 + \cdots,\\
j_3^*(\tau) &=& \biggl( \frac{\eta(\tau)}{\eta(3\tau)} \biggr)^{12} + 12 + 3^6 \biggl( \frac{\eta(3\tau)}{\eta(\tau)} \biggr)^{12} = \frac{1}{q} + 783q + 8672q^2 + 65367q^3 + \cdots,\\
j_5(\tau) &=& \biggl( \frac{\eta(\tau)}{\eta(5\tau)} \biggr)^6 + 6 = \frac{1}{q} + 9q + 10q^2 - 30q^3 + \cdots,\\
j_5^*(\tau) &=& \biggl( \frac{\eta(\tau)}{\eta(5\tau)} \biggr)^6 + 6 + 5^3 \biggl( \frac{\eta(5\tau)}{\eta(\tau)} \biggr)^6 = \frac{1}{q} + 134q + 760q^2 + 3345q^3 + \cdots.
\end{eqnarray*}

For $p$ = 2, 3, and 5, let $d$ be a positive integer such that $-d$ is congruent to a square modulo 4$p$, and $\mathcal{Q}_{d,p}$ the set of positive definite binary quadratic forms $Q(X,Y) = [a,b,c] = a X^2 + b X Y + c Y^2\ (a,b,c \in \mathbb{Z})$ of discriminant $-d$ with $a \equiv 0$ (mod $p$). Moreover, we fix an integer $\beta$ (mod $2p$) with $\beta^2 \equiv -d$ (mod $4p$) and denote by $\mathcal{Q}_{d,p,\beta}$ the set of quadratic forms $[a,b,c] \in \mathcal{Q}_{d,p}$ such that $b \equiv \beta$ (mod $2p$). For every positive integer $m$, let $\varphi_m(j_p^*)$ be a unique polynomial of $j_p^*$ satisfying $\varphi_m(j_p^*(\tau)) = q^{-m} + O(q)$. We define two modular trace functions:
\begin{eqnarray*}
\textbf{t}_m^{(p)}(d) &:=& \sum_{Q \in \mathcal{Q}_{d,p,\beta} / \Gamma_0(p)} \frac{1}{|\Gamma_0(p)_Q|} \varphi_m(j_p^*(\alpha_Q)),\\
\textbf{t}_m^{(p*)}(d) &:=& \sum_{Q \in \mathcal{Q}_{d,p} / \Gamma_0^*(p)} \frac{1}{|\Gamma_0^*(p)_Q|} \varphi_m(j_p^*(\alpha_Q)),
\end{eqnarray*}
where $\alpha_Q$ is the root of $Q(X,1) = 0$ in $\mathfrak{H}$. The definition of $\textbf{t}_m^{(p)}(d)$ is independent of $\beta$. In addition, we set $\textbf{t}_2^{(2*)}(0) := 5, \, \textbf{t}_2^{(3*)}(0) = \textbf{t}_2^{(5*)}(0) := 3, \, \textbf{t}_2^{(p*)}(-1) := -1, \, \textbf{t}_2^{(p*)}(-4) := -2, \, \textbf{t}_2^{(p*)}(d) := 0$ for $d < -4$ or $-d \not\equiv$ square (mod $4p$) ($p$ = 2, 3, 5). For the relation between two modular trace functions, see $\cite{Kim08}$.

\begin{rmk}
For $p = 1$, we put $j_1^*(\tau) := j(\tau) - 744 = \{(\eta(\tau)/\eta(2\tau))^8 + 2^8 (\eta(2\tau)/\eta(\tau))^{16}\}^3-744$ and $\textbf{t}_m(d) := \textbf{t}_m^{(1*)}(d)$.
\end{rmk}

\section{Proof of Theorem \ref{main1}}

We give a proof only for the case $p = 3$; the other cases are proved in the same way. 

\begin{dfn}
For every positive integer $t$, we define the operator $U_t$ by
\begin{eqnarray*}
\biggl(\sum a_n q^n\biggr) \biggr|U_t := \sum a_{tn} q^n.
\end{eqnarray*}
\end{dfn}
Then $U_t$ sends a modular form to a modular form of the same weight but raises the level in general. To prove Theorem \ref{main1}, we need the following theorem, which is a special case  $f = \varphi_m(j_p^*(\tau))$ of Theorem 1.1 in \cite{BF06}.

\begin{thm}
 \label{BF}
The function
\begin{eqnarray*}
g_m^{(p*)}(\tau) := \sum_{d > 0} \textbf{t}_m^{(p*)}(d) q^d + (\sigma_1(m) + p \sigma_1(m/p)) - \sum_{k |m}k q^{-k^2}
\end{eqnarray*}
$($where $\sigma_1(x) = 0$ if $x \not\in \mathbb{Z})$ is a meromorphic modular form of weight $3/2$, holomorphic outside the cusps, with respect to $\Gamma_0(4p)$, that is,
\begin{eqnarray*}
g_m^{(p*)}(\tau) \in M_{3/2}^{mer}(\Gamma_0(4p)).
\end{eqnarray*}
Here $M_k^{mer}(\Gamma)$ denotes the space of meromorphic modular forms of weight $k$ with respect to $\Gamma$.
\end{thm}

We prove Theorem\ref{main1}. For the modular form $f(\tau) = \sum a_n q^n$, we define the functions $\tilde{f}_0, \ \tilde{f}_1 \ and\  \tilde{f}_2$ by
\begin{eqnarray*}
\tilde{f}_0(\tau) &:=& \frac{1}{3} \left\{f(\tau) + f(\tau + \frac{1}{3}) + f(\tau + \frac{2}{3})\right\},\\
\tilde{f}_1(\tau) &:=& \frac{1}{3} \left\{f(\tau) + \zeta^{-1}f(\tau + \frac{1}{3}) + \zeta f(\tau + \frac{2}{3})\right\},\\
\tilde{f}_2(\tau) &:=& \frac{1}{3} \left\{f(\tau) + \zeta f(\tau + \frac{1}{3}) + \zeta^{-1}f(\tau + \frac{2}{3})\right\}
\end{eqnarray*}
where $\zeta = e^{2\pi i/3}$. For each $k \pmod{3}$, then $\tilde{f}_k$ has a Fourier expansion of the form $\tilde{f}_k(\tau) = \sum_{n \equiv k (3)}a_n q^n$, and it is also a modular form of the same weight. By Theorem \ref{BF},  we have
\begin{eqnarray*}
g_2^{(3*)}(\tau) = \sum_{d = -4}^{\infty} \textbf{t}_2^{(3*)}(d) q^d \in M_{3/2}^{mer}(\Gamma_0(12)).
\end{eqnarray*}
Now consider the modular form $g_2^{(3*)}(\tau) \cdot \theta_0(\tau)$ where $\theta_0(\tau) := \sum_{n \in \mathbb{Z}} q^{n^2} \in M_{1/2}(\Gamma_0(4))$. This form is of weight 2 and we have
\begin{eqnarray*}
g_2^{(3*)}(\tau) \cdot \theta_0(\tau) = \biggl(\sum_{d = -4}^{\infty}\textbf{t}_2^{(3*)}(d) q^d \biggr) \cdot \biggl(\sum_{r \in \mathbb{Z}}q^{r^2} \biggr) = \sum_{n = -4}^{\infty} \biggl(\sum_{r \in \mathbb{Z}} \textbf{t}_2^{(3*)}(n - r^2) \biggr) q^n \in M_2^{mer}(\Gamma_0(12)).
\end{eqnarray*}
Similarly, the product $g_2^{(3*)}(\tau) \cdot \theta_0(9\tau)$ is also a modular form of weight 2 and its Fourier expansion is
\begin{eqnarray*}
g_2^{(3*)}(\tau) \cdot \theta_0(9\tau) = \sum_{n = -4}^{\infty} \biggl(\sum_{r \in \mathbb{Z}} \textbf{t}_2^{(3*)}(n - (3r)^2) \biggr) q^n = \sum_{n = -4}^{\infty} \biggl(\sum_{r \equiv 0 (3)} \textbf{t}_2^{(3*)}(n - r^2) \biggr) q^n \in M_2^{mer}(\Gamma_0(36)).
\end{eqnarray*}
We put
\begin{eqnarray*}
F(\tau) &:=& \biggl(g_2^{(3*)}(\tau) \cdot \theta_0(\tau)\biggr) \bigg|U_4 =  \sum_{n = -1}^{\infty} \biggl(\sum_{r \in \mathbb{Z}} \textbf{t}_2^{(3*)}(4n - r^2) \biggr) q^n \in M_2^{mer}(\Gamma(12)), \\
G(\tau) &:=& \biggl(g_2^{(3*)}(\tau) \cdot \theta_0(9\tau)\biggr) \bigg|U_4 =  \sum_{n = -1}^{\infty} \biggl(\sum_{r \equiv 0 (3)} \textbf{t}_2^{(3*)}(4n - r^2) \biggr) q^n \in M_2^{mer}(\Gamma(36)). 
\end{eqnarray*}
Then $F(\tau)$ and $G(\tau)$ are meromorphic modular forms of weight 2. Moreover, for
\begin{eqnarray*}
j'_3(\tau) &=& \sum_{n = -1}^{\infty} n c_n^{(3)} q^n \in M_2^{mer}(\Gamma_0(3)),\\
E_2^{(3)}(\tau) &:=& \frac{1}{2}(3E_2(3\tau) - E_2(\tau)) = 1 + 12\sum_{n = 1}^{\infty} \sigma_1^{(3)}(n) q^n \in M_2(\Gamma_0(3)),
\end{eqnarray*}
(where the prime denotes $(2\pi i)^{-1} d/d\tau$ and $E_2(\tau) := 1 - 24\sum_{n = 1}^{\infty} \sigma_1(n) q^n$ is the Eisenstein series of weight 2), we put
\begin{eqnarray*}
H(\tau) := j'_3(\tau) - \frac{3}{2} E_2^{(3)}(\tau) = -\frac{1}{q} - \frac{3}{2} + \sum_{n = 1}^{\infty} (n c_n^{(3)} - 18\sigma_1^{(3)}(n)) q^n \in M_2^{mer}(\Gamma_0(3)).
\end{eqnarray*}
Then, the theorem in the case of $p = 3$ is equivalent to the following identities of modular forms:
\begin{eqnarray*}
2\tilde{H}_0(\tau) = -\tilde{G}_0(\tau), \ \ 2\tilde{H}_1(\tau) = \tilde{F}_1(\tau), \ \ 2\tilde{H}_2(\tau) = \tilde{F}_2(\tau).
\end{eqnarray*}
Since these modular forms are of weight 2 on $\Gamma(36)$, we see that, by the Riemann-Roch theorem, it is enough to check the coincidence of Fourier coefficients on both sides of the equalities up to $q^{3960}$. We checked this by using Mathematica and Pari-GP.\\

Similarly, we can show the equation $j_3^*(\tau) = j_3(\tau) - 3(j_3|U_3)(\tau)$, and we obtain $c_n^{(3*)} = c_n^{(3)} - 3c_{3n}^{(3)}$.

\section{Proof of Theorem \ref{main2}}

In this section, we give an overview of a proof. Since we can prove any case in the same way as \cite{Sam15}, we give a proof only for the case $p =$ 3. First, we prepare for a proof.

\begin{dfn}
The binary quadratic forms
\begin{eqnarray*}
 \left\{ \begin{array}{ll}
 [3, 0, d/12] & (-d \equiv 0 \pmod{12}), \\
 \lbrack 3, 1, (d+1)/12 \rbrack \ ,\  \lbrack 3, -1, (d+1)/12 \rbrack & (-d \equiv 1 \pmod{12}), \\ 
\lbrack 3, 2, (d+4)/12 \rbrack \ ,\  \lbrack 3, -2, (d+4)/12 \rbrack & (-d \equiv 4 \pmod{12}), \\
\lbrack 3, 3, (d+9)/12 \rbrack & (-d \equiv 9 \pmod{12}) \\ 
 \end{array} \right.
\end{eqnarray*}
are forms with discriminant $-d$ and are called the principal form of discriminant $-d$.
\end{dfn}

\begin{lem}
The following conditions are equivalent for a form $Q \in \mathcal{Q}_{d, 3}$:\\
$(1)$ There are $x, y \in \mathbb{Z}$ such that $Q(x, y) = 3$.\\
$(2)$ $Q$ is $\Gamma_0^*(3)$-equivalent to $[3, B, C]$ for some $B, C \in \mathbb{Z}$.\\
$(3)$ $Q$ is $\Gamma_0^*(3)$-equivalent to a principal form of discriminant $-d$.
\end{lem}

This lemma can be proved in the same way as Lemma 2.2 in \cite{Sam15}. The key theorem for the proof of Theorem\ref{main2} is the following.

\begin{thm}
 $($Laplace's method$)$. Suppose that $h(t)$ is a real-valued $C^2$-function defined on the interval $(a, b)$ $($with $a, b \in \mathbb{R}$$)$. If we further suppose that $h$ has a unique maximum at $t = c$ with $a < c < b$ so that $h'(c) = 0$ and $h''(c) < 0$, then, we have
\begin{eqnarray*}
\int_a^b e^{\lambda h(t)} dt \sim e^{\lambda h(c)} \biggl(\frac{-2\pi}{\lambda h''(c)} \biggr)^{1/2}
\end{eqnarray*}
as $\lambda \to \infty$.
\end{thm}

We prove Theorem \ref{main2}. By definition,
\begin{eqnarray*}
\textbf{t}_2^{(3*)}(d) &:=& \sum_{Q \in \mathcal{Q}_{d,3} / \Gamma_0^*(3)} \frac{1}{|\Gamma_0^*(3)_Q|} \varphi_2(j_3^*(\alpha_Q)).
\end{eqnarray*}
If $Q = [a, b, c]$ is the element of $\mathcal{Q}_{d, 3}$, we have
\begin{eqnarray*}
e^{2\pi i \alpha_Q} = \exp \biggl(2\pi i \biggl(\frac{-b + i \sqrt{d}}{2a} \biggr) \biggr) = \exp \biggl(-\frac{\pi i b}{a} \biggr) \exp \biggl(-\frac{\pi \sqrt{d}}{a} \biggr)
\end{eqnarray*}
and consequently;
\begin{eqnarray*}
\varphi_2(j_3^*(\alpha_Q)) &=& q^{-2} + O(q)\\
&=& \exp \biggl(\frac{2\pi i b}{a} \biggr) \exp \biggl(\frac{2\pi \sqrt{d}}{a} \biggr) + O\biggl(\exp \biggl(-\frac{\pi \sqrt{d}}{a} \biggr)\biggr).
\end{eqnarray*}
By this calculation, the contribution to $\textbf{t}_2^{(3*)}(d)$ comes only from classes of forms with $a = 3$. By Lemma 4.2, any such form is equivalent to a principal form, so that we have 
\begin{eqnarray*}
\textbf{t}_2^{(3*)}(d) = O\biggl(\exp \biggl(-\frac{\pi \sqrt{d}}{3}\biggr)\biggr) +  \exp \biggl(\frac{2\pi \sqrt{d}}{3} \biggr) \times \left\{ \begin{array}{ll} 1 & (d \equiv 0, 3\ \bmod 12), \\
 -1 & (d \equiv 8, 11 \bmod 12). \\
 \end{array} \right.
\end{eqnarray*}
Combining this formula with Theorem\ref{main1}, we obtain
\begin{eqnarray*}
c_n^{(3)} \sim \frac{1}{2n} \times \left\{ \begin{array}{ll}
 -\sum_{\substack{r \equiv 0 (3) \\ 4n \geq r^2}} \exp \bigl(2\pi \sqrt{4n - r^2}/3 \bigr) & (n \equiv 0 \bmod 3), \\
 \sum_{\substack{r \equiv 1, 2 (3) \\ 4n \geq r^2}} \exp \bigl(2\pi \sqrt{4n - r^2}/3 \bigr) & (n \equiv 1 \bmod 3), \\
 -\sum_{\substack{r \equiv 0 (3) \\ 4n \geq r^2}} \exp \bigl(2\pi \sqrt{4n - r^2}/3 \bigr) & (n \equiv 2 \bmod 3). \\
 \end{array} \right.
\end{eqnarray*}
For each $k = 0, 1, 2$, we consider the sum
\begin{eqnarray*}
S_n^{(k)} := \frac{3}{2\sqrt{n}} \sum_{\substack{r \equiv k (3) \\ 4n \geq r^2}} e^{\frac{4}{3}\pi \sqrt{n} \sqrt{1 - \frac{r^2}{4n}}} 
= \frac{3}{2\sqrt{n}} \sum_{\substack{l \in \mathbb{Z} \\ 4n \geq (3l + k)^2}} e^{\frac{4}{3}\pi \sqrt{n} \sqrt{1 - \frac{(3l + k)^2}{4n}}},
\end{eqnarray*}
and view this sum as a Riemann sum for the function $t \mapsto e^{4\pi \sqrt{n} \sqrt{1 - t^2}/3}$ : $(-1, 1) \to \mathbb{R}$. We can show that $S_n^{(k)}$ is asymptotic to the corresponding Riemann integral $J_n$ where
\begin{eqnarray*}
J_n := \int_{-1}^1 e^{4\pi \sqrt{n} \sqrt{1 - t^2}/3} dt.
\end{eqnarray*}
(For further detail, see \cite{Sam15}). Moreover, applying Laplace's method to the case $\lambda = \sqrt{n}$ and $h(t) = 4\pi \sqrt{1 - t^2}/3$ on $(-1, 1)$, we have
\begin{eqnarray*}
J_n \sim e^{\sqrt{n} \cdot 4 \pi/3} \cdot \biggl(\frac{-2\pi}{-4\pi \sqrt{n}/3} \biggr)^{1/2} = \frac{\sqrt{3}}{\sqrt{2} n^{1/4}} e^{4\pi \sqrt{n}/3}.
\end{eqnarray*}
Putting these asymptotic formulas together, we obtain
\begin{eqnarray*}
c_n^{(3)} &\sim& \frac{1}{3 \sqrt{n}} \times \left\{ \begin{array}{ll}
 -S_n^{(0)} & (n \equiv 0 \bmod 3), \\
 S_n^{(1)} + S_n^{(2)} & (n \equiv 1 \bmod 3), \\
 -S_n^{(0)} & (n \equiv 2 \bmod 3), \\
\end{array} \right. \\
&\sim& \frac{e^{4\pi \sqrt{n}/3}}{\sqrt{6} n^{3/4}} \times \left\{ \begin{array}{ll}
 -1 & (n \equiv 0 \bmod 3), \\
 2 & (n \equiv 1 \bmod3), \\
 -1 & (n \equiv 2 \bmod3) \\
\end{array} \right. 
\end{eqnarray*} 
as $n \to \infty$.

\newpage

\section{Tables of $\textbf{t}_m^{(p*)}(d)$ and $\textbf{t}_m^{(p)}(d)$ \ $(-4 \leq d \leq 50)$}

\begin{table}[h]
\begin{minipage}{9cm}
  \begin{tabular}[t]{|c|r|r|r|r|} \hline
    $d$ & $\textbf{t}_1^{(2*)}(d)$ & $\textbf{t}_2^{(2*)}(d)$ & $\textbf{t}_1^{(2)}(d)$ & $\textbf{t}_2^{(2)}(d)$\\ \hline \hline
    $-$4 & 0 & $-$2 & 0 & $-$4  \\ \hline
    $-$1 & $-$1 & $-$1 & $-$1 & $-$1 \\ \hline
    0 & 1 & 5 & 2 & 10 \\ \hline
    4 & $-$26 & 518 & $-$52 & 1036 \\ \hline
    7 & $-$23 & $-$8215 & $-$23 & $-$8215 \\ \hline
    8 & 76 & 7180 & 152 &14360 \\ \hline
    12 & $-$248 & 52760 & $-$496 & 105520 \\ \hline
    15 & $-$1 & $-$385025 & $-$1 & $-$385025 \\ \hline
    16 & 518 & 287710 & 1036 & 575420 \\ \hline
    20 & $-$1128 &1263640 & $-$2256 & 2527280 \\ \hline
    23 & $-$94 & $-$6987870 & $-$94 & $-$6987870 \\ \hline
    24 & 2200 & 4831256 & 4400 & 9662512 \\ \hline
    28 & $-$4096 & 16572370 & $-$8192 & 33144740 \\ \hline
    31 & 93 & $-$78987171 & 93 & $-$78987171 \\ \hline
    32 & 7180 & 52263100 & 14360 & 104526200 \\ \hline
    36 & $-$12418 & 153553438 & $-$24836 & 307106876 \\ \hline
    39 & $-$236 & $-$663068908 & $-$236 & $-$663068908 \\ \hline
    40 & 20632 & 425670680 & 41264 & 851341360 \\ \hline
    44 & $-$33512 & 1122593352 & $-$67024 & 2245186704 \\ \hline
    47 & 235 & $-$4515675925 & 235 & $-$4515675925 \\ \hline
    48 & 53256 & 2835914280 & 106512 & 5671828560 \\ \hline
\end{tabular}
\end{minipage}
\hfill
\begin{minipage}{9cm}
  \begin{tabular}[t]{|c|r|r|r|r|} \hline
    $d$ & $\textbf{t}_1^{(3*)}(d)$ & $\textbf{t}_2^{(3*)}(d)$ & $\textbf{t}_1^{(3)}(d)$ & $\textbf{t}_2^{(3)}(d)$\\ \hline \hline
    $-$4 & 0 & $-$2 & 0 & $-$2  \\ \hline
    $-$1 & $-$1 & $-$1 & $-$1 & $-$1 \\ \hline
    0 & 1 & 3 & 2 & 6 \\ \hline
    3 & $-$7 & 33 & $-$14 & 66 \\ \hline
    8 & $-$34 & $-$410 & $-$34 & $-$410 \\ \hline
    11 & 22 & $-$1082 & 22 &$-$1082 \\ \hline
    12 & 26 & 1428 & 52 & 2856 \\ \hline
    15 & $-$69 & 3195 & $-$138 & 6390 \\ \hline
    20 & $-$116 & $-$11892 & $-$116 & $-$11892 \\ \hline
    23 & 115 & $-$22797 & 115 & $-$22797 \\ \hline
    24 & 174 & 28710 & 348 & 57420 \\ \hline
    27 & $-$241 & 53223 & $-$482 & 106446 \\ \hline
    32 & $-$410 & $-$140222 & $-$410 & $-$140222 \\ \hline
    35 & 492 & $-$240500 & 492 & $-$240500 \\ \hline
    36 & 492 & 287244 & 984 & 574488 \\ \hline
    39 & $-$705 & 477567 & $-$1410 & 955134 \\ \hline
    44 & $-$1060 & $-$1081096 & $-$1060 & $-$1081096 \\ \hline
    47 & 1272 & $-$1718792 & 1272 & $-$1718792 \\ \hline
    48 & 1442 & 2004918 & 2884 & 4009836 \\ \hline
\end{tabular}
\end{minipage}
\end{table}

\begin{table}[h]
  \begin{tabular}{|c|r|r|r|r|} \hline
    $d$ & $\textbf{t}_1^{(5*)}(d)$ & $\textbf{t}_2^{(5*)}(d)$ & $\textbf{t}_1^{(5)}(d)$ & $\textbf{t}_2^{(5)}(d)$\\ \hline \hline
    $-$4 & 0 & $-$2 & 0 & $-$2  \\ \hline
    $-$1 & $-$1 & $-$1 & $-$1 & $-$1 \\ \hline
    0 & 1 & 3 & 2 & 6 \\ \hline
    4 & $-$8 & $-$6 & $-$8 & $-$6 \\ \hline
    11 & $-$12 & $-$124 & $-$12 & $-$124 \\ \hline
    15 & $-$19 & 93 & $-$38 & 186 \\ \hline
    16 & $-$6 & $-$270 & $-$6 & $-$270 \\ \hline
    19 & 20 & 132 & 20 & 132 \\ \hline
    20 & 6 & 268 & 12 & 536 \\ \hline
    24 & $-$44 & 216 & $-$44 & 216 \\ \hline
    31 & $-$39 & $-$1863 & $-$39 & $-$1863 \\ \hline
    35 & $-$44 & 1668 & $-$88 & 3336 \\ \hline
    36 & 20 & $-$3054 & 20 & $-$3054 \\ \hline
    39 & 53 & 1653 & 53 & 1653 \\ \hline
    40 & 56 & 2868 & 112 & 5736 \\ \hline
    44 & $-$136 & 2416 & $-$136 & 2416 \\ \hline
\end{tabular}
\end{table}

\newpage

\begin{ack}
The authors are grateful to Professor Masanobu Kaneko for his helpful comments and encouraging our studies. 
\end{ack}
%参考文献

\noindent
T. Matsusaka: Graduate School of Mathematics, Kyushu University, Motooka 744, Nishi-ku Fukuoka 819-0395, Japan\\
e-mail: toshikimatsusaka@gmail.com\\

\noindent
R. Osanai:\\
e-mail: ryotaroosanai@gmail.com


\begin{thebibliography}{99}
\bibitem{BF06} J. H. Bruinier and J Funke,
Traces of CM values of modular functions. J. Reine Angew. Math. \textbf{594}, (2006), 1--33.
\bibitem{BO13} J. H. Bruinier and K. Ono, 
Algebraic formulas for the coefficients of half-integral weight harmonic weak Maass forms. Adv. Math. \textbf{246}, (2013), 198--219.
\bibitem{Kim08} D. Choi, D. Jeon, S. Y. Kang, and C. H. Kim,
Exact formulas for traces of singular moduli of higher level modular functions. J. Number Theory \textbf{128}, (3), (2008), 700--707.
\bibitem{DM13} M. Dewar and M. R. Murty, 
A derivation of the Hardy-Ramanujan formula from an arithmetic formula. Proc. Amer. Math. Soc. \textbf{141}, (2013), no. 6, 1903--1911.
\bibitem{DM132} M. Dewar and M. R. Murty, 
An asymptotic formula for the coefficients of $j(\tau)$. Int. J. Number Theory \textbf{9}, (2013), no. 3, 641--652. 
\bibitem{HR18} G. H. Hardy and S. Ramanujan, 
Asymptotic formulae in combinatory analysis. Proc. London Math. Soc. (2) \textbf{17}, (1918), 75--115.
\bibitem{Kan95} M. Kaneko,
The Fourier coefficients and the singular moduli of the elliptic modular function $j(\tau)$. Mem. Fac. Engrg. Design Kyoto Inst. Tech. Ser. Sci. Tech. \textbf{44}, (1996), 1--5.
\bibitem{Sam15} M. R. Murty and K. Sampath, 
On the asymptotic formula for the Fourier coefficients of $j$-function. Kyushu J. Math. \textbf{70}, (2016).
\bibitem{Ohta09} K. Ohta,
Formulas for the Fourier coefficients of some genus zero modular functions. Kyushu J. Math. \textbf{63}, (2009), 1--15.
\bibitem{Pet32} H. Petersson,
\"Uber die Entwicklungskoeffizienten der automorphen Formen. Acta Math. \textbf{58}, (1), (1932), 169--215.
\bibitem{Rad} H. Rademacher,
The Fourier coefficients of the modular invariant $J(\tau)$. Amer. J. Math.,
\textbf{60}, (2), (1938), 501--512.
\bibitem{Zag02} D. Zagier,
Traces of singular moduli. Motives, Polylogarithms and Hodge theory, Part I (Irvine, CA, 1998), 211--244, Int. Press Lect. Ser. 3, I, 2002.\\
\end{thebibliography}
\end{document}